\documentclass{ifacconf}

\usepackage{graphicx}      
\usepackage{natbib}        

\usepackage{amsmath}
\usepackage{amssymb}

\newcommand{\XX}{{\mathcal X }}
\newcommand{\YY}{{\mathcal Y }}

\newcommand{\UU}{{\mathcal U}}
\newcommand{\LL}{{\mathcal L}}
\newcommand{\II}{{\mathcal I}}

\newcommand{\abs}[1]{\lvert#1\rvert}

\newcommand{\dotex}{{\frac{d}{dt}}}

\newcommand{\eps}{\varepsilon}
\DeclareMathOperator{\sign}{sign}

\begin{document}
\begin{frontmatter}

\title{A Separation Principle on Lie Groups}


\author[First]{S. Bonnabel}
\author[First]{P. Martin}
\author[First]{P. Rouchon}
\author[Second]{E. Sala\"un}

\address[First]{Math\'ematiques et syst\`emes, Mines ParisTech, 75272 Paris Cedex~06, France (e-mail: [silvere.bonnabel,philippe.martin,pierre.rouchon]@mines-paristech.fr)}
\address[Second]{School of Aerospace Engineering, Georgia Institute of Technology, Atlanta, GA 30332-0150, USA (e-mail:erwan.salaun@gatech.edu)}

\begin{abstract}                
For linear time-invariant systems, a separation principle holds:  stable observer and stable state feedback can be designed for the time-invariant system, and the combined observer and feedback will be stable. For non-linear systems, a local separation principle holds around steady-states, as the linearized system is time-invariant. 
This paper addresses the issue of a non-linear separation principle on Lie groups.  For invariant systems on Lie groups, we prove there exists a large set of (time-varying) trajectories around which the linearized observer-controler system is time-invariant, as soon as a symmetry-preserving observer is used. Thus a separation principle holds around those trajectories.  The theory is illustrated by a mobile robot example, and the developed ideas are then  extended to a class of Lagrangian mechanical systems on Lie groups  described by Euler-Poincar\'e equations.\end{abstract}


\end{frontmatter}
\section{Introduction}

The celebrated separation principle plays a key role in the linear theory of control. It states that the problem of designing an optimal stable feedback controller can be broken into two parts: designing an optimal stable observer meant to feed an optimal  stable controller.  As a particular application of this principle it is proved that for a  \emph{linear time-invariant} system, combining  a stable observer and a stable controller yields  a stable closed-loop system.

When the system is non-linear, the separation principle does not holds. However a local separation principle can always be stated around steady-states, as the linearized system around steady-states is linear and time-invariant. But in general it does not hold around other types of trajectories. In this paper, we consider invariant systems on Lie groups, and we state a local separation principle around a large class of trajectories that are not necessarily steady-states.  Note that  separation principles  for special classes of non-linear systems have been addressed in e.g.  \cite{Khalil-99,Gauthier-Kupka-92}, \cite{maithripala2005}.

 Various systems of
engineering interest can be  modeled as invariant systems on Lie groups, mainly
 cart-like  vehicles and mechanical systems such as rigid bodies in space.  There is an extensive literature on control on Lie groups (see e.g.  \cite{bullo-murray-auto99,grizzle-marcus-ieee85,respondek-tall-scl02,morin-ieee03} and \cite{Sussmann-72} as one of the pioneering papers), but general methods for observer design on Lie groups have only been introduced recently (\cite{arxiv-08,lageman10,mahony-et-al-IEEE}). In this paper, a link is established between those two fields or reasearch and a separation principle on Lie groups is stated.

In Section 2, we recall the linear separation principle, and the local non-linear separation principle around steady states. In Section 3, we consider a  left-invariant system on a Lie group and we build a symmetry-preserving observer. It has been proved that, for such observers there exists a large set of trajectories around which the estimation  error is time-invariant: the so-called permanent trajectories (\cite{arxiv-07,arxiv-08}). In this paper, we prove that the tracking error is also time-invariant when the reference trajectory is a permanent trajectory. As a result, the \emph{linearized} closed-loop observer-controler system is \emph{time-invariant} around permanent trajectories. This implies a (local) separation principle around permanent trajectories. This result advocates that the recently introduced permanent trajectories for systems on Lie groups  generalize steady states for general non-linear systems. Indeed, if the observer is not a symmetry-preserving observer, the linearized closed-loop observer-controler system around permanent  trajectories is time-varying (in general), implying that  a local separation principle does not hold.

In Section 4, we consider the well-known problem of localizing from landmarks using sonar (see e.g. \cite{betke-et-al-02}). Observability follows from triangulation, and position estimation is generally achieved via the widespread
 Extended Kalman filter (see e.g. \cite{Roumeliotis}). For this system we derive a class of symmetry-preserving observers that converge around \emph{any} permanent trajectories (i.e. circles and lines with constant speed). The separation principle derived in this paper allows to prove that, as soon as the control scheme is stable around those trajectories,  the
 closed-loop system is also stable around those trajectories and the eigenvalues of the closed-loop linearized system are those of the observer together with those of the controller. It is interesting to note that this stability result does not hold for an Extended Kalman Filter. Note that symmetry-preserving observers have already been used for localization from landmarks in \cite{vasconcelos}. 

Besides the two main contribution of this paper, which are to derive a (local) separation principle for  invariant systems on Lie groups, and to derive an observer-controler for the localization problem from landmarks  with guaranteed convergence properties around a large set of trajectories, we consider in Section 5 a (very) particular observation problem for simple mechanical systems on Lie groups whose motion is described by Euler-Poincar\'e equations. Such systems have been extensively studied in the literature (see e.g. \cite{bullo-lewis-article00,bullo-murray-auto99,marsden-ratiu-book94}). By a similar token as the one of Section 2 we prove a local separation principle around permanent trajectories.

\section{The linear separation principle}\label{sec:lsprinciple}
Consider the system
\begin{align}
\label{eq:x}\dot x &=f(x,u)\\
\label{eq:y}y &=h(x,u),
\end{align}
where $(x,u,y)$ belongs to an open subset $\XX\times\UU\times\YY\subset\Rset^n\times\Rset^m\times\Rset^p$. We would like to track the reference trajectory
\begin{align}
\label{eq:xr}\dot x_r &=f(x_r,u_r)\\
\label{eq:yr}y_r &=h(x_r,u_r)
\end{align}
using only the measured output~$y$. In other words we want to stabilize the equilibrium point $(\bar\eta_x,\bar\eta_u):=(0,0)$ of the error system
\begin{align}
\label{eq:etax}\dot\eta_x &=f(x_r+\eta_x,u_r+\eta_u)-f(x_r,u_r)\\
\label{eq:etay}\eta_y &=h(x_r+\eta_x,u_r+\eta_u)-h(x_r,u_r),
\end{align}
where $\eta_x:=x-x_r$, $\eta_u:=u-u_r$ and $\eta_y:=y-y_r$. Here we are interested only in local stability, i.e. we want to stabilize the linearized error system
\begin{align}
\label{eq:xix}\dot\xi_x &=\partial_1f(x_r,u_r)\xi_x+\partial_2f(x_r,u_r)\xi_u &&=A\xi_x+B\xi_u\\
\label{eq:xiy}\xi_y &=\partial_1h(x_r)\xi_x+\partial_2h(x_r,u_r)\xi_u &&=C\xi_x+D\xi_u.
\end{align}
The notation $\partial_i$ stands for the derivative with respect to the $i^{th}$ argument.

Notice the system~\eqref{eq:etax}-\eqref{eq:etay} as well as the matrices $A,B,C,D$ are in general not time-invariant unless the reference trajectory $(x_r,u_r) $ is an equilibrium point, i.e. is constant and such that $f(x_r,u_r)=0$.

\subsection{Linear controller with linear observer}\label{sec:linobs}
To stabilize the linearized system, one can use the linear controller-observer
\begin{align}
\label{eq:lincon}\xi_u &=-K\hat\xi_x\\
\nonumber\dot{\hat\xi}_x &=A\hat\xi_x+B\xi_u-L(C\hat\xi_x+D\hat\xi_u-\xi_y)\\
\label{eq:linobs}&=(A-BK-LC+BLD)\hat\xi_x+L\xi_y,
\end{align}
where the $m\times n$ matrix~$K$ and $n\times p$ matrix~$L$ are to be chosen.
Indeed, setting $e_x:=\hat\xi_x-\xi_x$ the closed-loop system
\begin{align}
\label{eq:sepeta}\dot\xi_x &=(A-BK)\xi_x-BKe_x\\
\label{eq:sepe}\dot e_x &=(A-LC)e_x
\end{align}
has a triangular structure hence its eigenvalues are those of $A-BK$ together with those of $A-LC$. If $(x_r,u_r) $ is an equilibrium point, it can easily be stabilized by placing separately the eigenvalues of $A-BK$ through~$K$ and those of $A-LC$ through~$L$ (provided of course $(A,B)$ is controllable and $(A,C)$ is observable). This is the well-known separation principle for linear time-invariant systems. The result still applies if the reference trajectory is slowly-varying, i.e. $f\bigl(x_r(t),u_r(t)\bigr)\approx0$ for all~$t$; $K$ and~$L$ may also depend on the reference trajectory to provide gain-scheduling.

This shows that the control law
\begin{align}
\label{eq:linconf}\eta_u &=-K\hat\xi_x\\
\label{eq:linobsf}\dot{\hat\xi}_x &=(A-BK-LC+BLD)\xi_x+L\eta_y,
\end{align}
where the input~$\xi_y$ and output~$\xi_u$ of~\eqref{eq:lincon}-\eqref{eq:linobs} are replaced by $\eta_y$ and~$\eta_u$, locally stabilizes the equilibrium point $(\bar\eta_x,\bar\eta_u):=(0,0)$ of~\eqref{eq:etax}. In other words
\begin{align*}
u &=u_r-K\hat\xi_x\\
\dot{\hat\xi}_x &=(A-BK-LC+BLD)\hat\xi_x+L(y-y_r)
\end{align*}
stabilizes \eqref{eq:x} around the reference trajectory~\eqref{eq:xr}-\eqref{eq:yr} using only the measured output~\eqref{eq:y}.

\subsection{Linear controller with extended observer}\label{sec:extobs}
Instead of the linear observer of the previous section, the so-called ``extended'' observer
\begin{align*}
\dot{\hat x} &=f(\hat x,u)-L\bigl(h(\hat x,u)-y\bigr)
\end{align*}
can be used. The $n\times p$ gain matrix $L$ may depend on~$\hat x$ to provide gain-scheduling. We assume this observer converges, i.e. the origin of the error system
\begin{align*}
\dot\eps_x &=f(\hat x,u)-L\bigl(h(\hat x,u)-h(\hat x-\eps_x,u)\bigr) - f(\hat x-\eps_x,u)
\end{align*}
is stable, where $\eps_x:=\hat x-x$ is the observation error. Once again we are interested only in local stability, i.e. we assume only the stability of the linearized error system
\begin{align}
\label{eq:ex}\dot e_x &=\bigl(\partial_1f(\hat x,u)-L\partial_1h(\hat x,u)\bigr)e_x.
\end{align}
Notice it is usually not easy to ensure even this local stability, unless the reference trajectory $(x_r,u_r)$ is an equilibrium point. In this case \eqref{eq:ex} can be further linearized around $(x_r,u_r)$, which yields~\eqref{eq:sepe}.

The control law
\begin{align*}
u &=u_r-K(\hat x-x_r)\\
\dot{\hat x} &=f(\hat x,u)-L\bigl(h(\hat x,u)-y\bigr)
\end{align*}
then locally stabilizes \eqref{eq:x} around the reference trajectory~\eqref{eq:xr}-\eqref{eq:yr} using only the measured output~\eqref{eq:y}. Indeed,
\begin{align*}
\eta_u=-K(\hat x-x+x-x_r)=-K(\eps_x+\eta_x),
\end{align*}
hence $\xi_u=-K(e_x+\xi_x)$, so that the linearized closed-loop system
\begin{align*}
\dot\xi_x &=(A-BK)\xi_x-BKe_x\\
\dot e_x &=\bigl(\partial_1f(\hat x,u)-L\partial_1h(\hat x,u)\bigr)e_x
\end{align*}
is clearly stable.

\section{A separation principle for invariant systems}
We again consider the system~\eqref{eq:x}-\eqref{eq:y}, but now $\XX$ is a Lie group of dimension~$n$ with identity~$e$ and group law
\[ x_2x_1:=\varphi(x_2,x_1). \]

\begin{defn}
Let $\Sigma$ be an open set (or more generally a manifold). A \emph{transformation group} on~$\Sigma$ is a smooth map
\[ (x,\xi)\in\XX\times\Sigma\mapsto\phi(x,\xi)\in\Sigma \]
such that
\begin{itemize}
 \item $\phi(e,\xi)=\xi$ for all~$\xi$
 \item $\phi\bigl(x_2,\phi(x_1,\xi)\bigr)=\phi(x_2x_1,\xi)$ for all $x_1,x_2,\xi$.
\end{itemize}
\end{defn}
By construction $\phi(x,\cdot)$ is a diffeomorphism on~$\Sigma$ for
all~$x$.

Consider then the transformation group on $\XX\times\UU\times\YY$ defined by
$\phi\bigl(x_0,(x,u,y)\bigr):=\bigl(\varphi(x_0,x),\psi(x_0,u),\varrho(x_0,y)\bigr)$. We will assume in the sequel that the system~\eqref{eq:x}-\eqref{eq:y} enjoys the following important invariance property.
\begin{defn}\label{def:invsys}
The system~\eqref{eq:x}-\eqref{eq:y} is \emph{invariant} by the transformation group if for all $x_0,x,u$
\begin{itemize}
 \item $f\bigl(\varphi(x_0,x),\psi(x_0,u)\bigr)=\partial_2\varphi(x_0,x)f(x,u)$
 \item $h\bigl(\varphi(x_0,x),\psi(x_0,u)\bigr)=\varrho\bigl(x_0,h(x,u)\bigr)$.
\end{itemize}
\end{defn}
With $(X,U,Y):=\bigl(\varphi(x_0,x),\psi(x_0,u),\varrho(x_0,u)\bigr)$ this reads
\begin{align*}
\dot X &=f(X,U)\\
Y &=h(X,U),
\end{align*}
i.e. \eqref{eq:x}-\eqref{eq:y} is left unchanged by the transformation group.

Finally define the following ``products'' which provide very compact notations: for $x_0\in\XX$, $(x,u,y)\in\XX\times\UU\times\YY$ and $\xi\in T_x\XX$ ($T_x\XX$ is the tangent space of $\XX$ at~$x$)
\begin{align}
x_0u &:=\psi(x_0,u)\\
x_0y &:=\varrho(x_0,y)\\
\label{eq:x0xi}x_0\xi &:=\partial_2\varphi(x_0,x)\xi\in T_{x_0x}\XX\\
\xi x_0 &:=\partial_1\varphi(x,x_0)\xi\in T_{xx_0}\XX.
\end{align}
Invariance then reads $f(x_0x,x_0u)=x_0f(x,u)$ and $h(x_0x,x_0u)=x_0h(x,u)$.

We would like the invariant system~\eqref{eq:x}-\eqref{eq:y} to track the reference trajectory~\eqref{eq:xr}-\eqref{eq:yr}. Instead of using as before the ``linear'' errors $x-x_r$, $u-u_r$ and~$y-y_r$, we consider the errors
$\eta_x:=x_r^{-1}x$, $\eta_u:=x_r^{-1}u-x_r^{-1}u_r$ and $\eta_y:=x_r^{-1}y-x_r^{-1}y_r$. These errors are invariant in the sense that $(x_0x_r)^{-1}(x_0x)=x_r^{-1}x$, $(x_0x_r)^{-1}(x_0u)=x_r^{-1}u$ and $(x_0x_r)^{-1}(x_0y)=x_r^{-1}y$.
The error system is then given by
\begin{align}
\nonumber\dot\eta_x &=x_r^{-1}f(x,u)-\bigl(x_r^{-1}f(x_r,u_r)x_r^{-1}\bigr)x\\
\nonumber&=\eta_xf(e,x^{-1}u)-f(e,x_r^{-1}u_r)\eta_x\\
\label{eq:etaxLie}&=\eta_xf(e,\eta_x^{-1}(I_r+\eta_u))-f(e,I_r)\eta_x\\
\nonumber\eta_y &=x_r^{-1}h(x,u)-x_r^{-1}h(x_r,x_r)\\
\label{eq:etayLie}&=h(\eta_x,I_r+\eta_u)-h(e,I_r),
\end{align}
and we want to stabilize its equilibrium point $(\bar\eta_x,\bar\eta_u):=(e,0)$. We have used the important invariant quantity
$$I_r:=x_r^{-1}u_r,$$
which leads to the following definition of a class of remarkable trajectories.
\begin{defn}
The trajectory $(x_r,u_r)$ is \emph{permanent} if $I_r$ is constant (i.e. independent of time).
\end{defn}

Once again we are interested only in local stability, i.e. we want to stabilize only the linearization of~\eqref{eq:etaxLie}-\eqref{eq:etayLie} around~$(\bar\eta_x,\bar\eta_u):=(e,0)$.
Recall that to linearize around~$e$ the map $\eta_x\in\XX\mapsto g(\eta_x)$, one can write $\eta_x=\exp(s\xi_x)$ where  $\exp$ is the exponential map of the Lie group~$\XX$, $s\in\Rset$ and~$\xi_x\in T_e\XX$ ($T_e\XX$ is identified with the Lie algebra of~$\XX$). The linearization of~$g$ is then given by
\begin{align*}
\left.\frac{d}{ds}g\bigl(\exp(s\xi_x)\bigr)\right|_{s=0}=\partial_1g(e)\xi_x.
\end{align*}
Applying this to~\eqref{eq:etaxLie}-\eqref{eq:etayLie}, we find the linearized error system is
\begin{align}
\nonumber\dot\xi_x &=\xi_xf(e,I_r)-f(e,I_r)\xi_x+\partial_2f(e,I_r)(\xi_u-\xi_xI_r)\\
\label{eq:xixLie}&=A\xi_x+B\xi_u\\
\nonumber\xi_y &=\partial_1h(e,I_r)\xi_x+\partial_2h(e,I_r)\xi_u\\
\label{eq:xiyLie}&=C\xi_x+D\xi_u,
\end{align}
where we have defined the ``products''
\begin{align*}
\xi u &:=\partial_1\psi(e,u)\xi\\
\xi\zeta &:=\bigl(\partial_1\partial_2\varphi(e,e)\xi\bigr)\zeta=\bigl(\partial_2\partial_1\varphi(e,e)\zeta\bigr)\xi
\end{align*}
for $u\in\UU$ and $\xi,\zeta\in T_e\XX$.

Notice the system~\eqref{eq:etaxLie}-\eqref{eq:etayLie} as well as the matrices $A,B,C,D$ are in general not time-invariant unless the reference trajectory $(x_r,u_r) $ is permanent. Since \eqref{eq:xixLie}-\eqref{eq:xiyLie} has the same form as~\eqref{eq:xix}-\eqref{eq:xiy}, we will be able to extend the usual separation principle around equilibrium points to a separation principle around permanent trajectories. The main benefit is that there are many more permanent trajectories thant equilibrium points.

\subsection{Linear controller with linear observer}
To stabilize the linearized error system~\eqref{eq:xixLie}-\eqref{eq:xiyLie} we can use a linear controller-observer of the form~\eqref{eq:lincon}-\eqref{eq:linobs}. As in section~\ref{sec:linobs}, the control law~\eqref{eq:linconf}-\eqref{eq:linobsf} with $\eta_x,\eta_u$ replacing $\xi_y,\xi_u$ stabilizes the equilibrium point $(\bar\eta_x,\bar\eta_u):=(e,0)$ of~\eqref{eq:etaxLie}. Using the definitions of $\eta_u,\eta_y$ and the fact that $x_r^{-1}(x_ru)=(x_r^{-1}x_r)u=u$ we eventually find that
\begin{align*}
u &=x_r(I_r-K\hat\xi_x)\\
\dot{\hat\xi}_x &=(A-BK-LC+BLD)\hat\xi_x+L\bigl(x_r^{-1}y-h(e,I_r)\bigr)
\end{align*}
stabilizes \eqref{eq:x} around the reference trajectory~\eqref{eq:xr}-\eqref{eq:yr} using only the measured output~\eqref{eq:y}.

\subsection{Linear controller with invariant observer}\label{sec:extobsLie}
Following~\cite{arxiv-07,arxiv-08}, we can easily design an invariant extended observer for invariant systems, using the
invariant output error $\eps_y:=\hat x^{-1}h(\hat x,u)-\hat x^{-1}y$ instead of the ``linear'' output error $h(\hat x,u)-y$.
Such an observer reads
\begin{align*}
\dot{\hat x} &=f(\hat x,u)-\hat xL\bigl(\hat x^{-1}h(\hat x,u)-\hat x^{-1}y\bigr),
\end{align*}
where the gain matrix $L$ may depend on $\hat I:=\hat x^{-1}u$ and on the output error~$\eps_y$.


The error system for the (invariant) state observation error $\eps_x:=x^{-1}\hat x$ is
\begin{align*}
\dot\eps_x &=x^{-1}\Bigl(f(\hat x,u)+\hat xL\bigl(\hat x^{-1}h(x,u)-\hat x^{-1}h(\hat x,u)\bigr)\Bigr)\\
&\quad-x^{-1}f(x,u)x^{-1}\hat x\\
&=\eps_xf(e,\hat I)-f(e,\eps_x\hat I)\eps_x-\eps_xL\bigl(h(\eps_x^{-1},\hat I)-h(e,\hat I)\bigr),
\end{align*}
and its linearization around the origin is
\begin{align*}
\dot e_x &=e_xf(e,\hat I)-f(e,\hat I)e_x-\partial_2f(e,\hat I)e_x-L\partial_1h(e,\hat I)e_x.
\end{align*}
It is then easy to tune the gain matrix~$L$ so that the observer converges at least around permanent trajectories.


Mimicking section~\ref{sec:extobs} the control law
\begin{align*}
u &=x_r\bigl(I_r-K\kappa(x_r^{-1}\hat x)\bigr)\\
\dot{\hat x} &=f(\hat x,u)+\hat xL\bigl(\hat x^{-1}y-h(e,\hat I)\bigr)
\end{align*}
then locally stabilizes \eqref{eq:x} around the reference trajectory~\eqref{eq:xr}-\eqref{eq:yr} using only the measured output~\eqref{eq:y}. Here $\kappa$ is any map such that
\begin{align*}
\left.\frac{d}{ds}\kappa\bigl(\exp(s\xi)\bigr)\right|_{s=0}=\xi,
\end{align*}
for instance a determination of the logarithm map of the group (i.e. the inverse of the exponential map of the group).
Indeed,
$x_r^{-1}u=I_r-K\kappa(x_r^{-1}xx^{-1}\hat x)$, i.e.
$\eta_u=-K\kappa(\eta_x\eps_x)$.
Linearizing, this yields
$\xi_u=-K(e_x+\xi_x)$, so that the linearized closed-loop system
\begin{align*}
\dot\xi_x &=(A-BK)\xi_x-BKe_x\\
\dot e_x &=e_xf(e,\hat I)-f(e,\hat I)e_x-\partial_2f(e,\hat I)e_x-L\partial_1h(e,\hat I)e_x
\end{align*}
is clearly stable.

\section{Example: output control of a wheeled robot from sonar landmarks}\label{sec:example}
Consider the simple planar non-holonomic model of a wheeled robot
\begin{align*}
\dot x &=uC\theta\\
\dot y &=uS\theta\\
\dot\theta &=uv,
\end{align*}
where the control inputs are the velocity~$u$ and $v$ the tangent of the steering angle; $C\theta$ and $S\theta$ stand for $\cos\theta$ and $\sin\theta$. We take for measurements the (square of the) distances to $p\geq3$ planar landmarks with known coordinates~$(x_i,y_i)$,
\begin{align*}
\lambda_i &=(x-x_i)^2+(y-y_i)^2,\qquad i=1,\ldots,p.
\end{align*}
The measurement system could be for instance a set of ultrasonic emitters located at~$(x_i,y_i)$, while the robot is equipped with ultrasonic receivers measuring the time of flight of the sonic waves.

Notice the only equilibrium points are given by $\bar u=0$ and $\bar x,\bar y,\bar\theta,\bar v$ constant, i.e. the robot is at rest. Moreover the linearized system around these points is clearly not controllable. Hence it is not possible to use the linear separation principle of section~\ref{sec:lsprinciple}.

It is easy to check the composition law
\begin{align*}
\begin{pmatrix}x_0\\y_0\\\theta_0\end{pmatrix}\begin{pmatrix}x\\y\\\theta\end{pmatrix}
&:=\begin{pmatrix}xC\theta_0-yS\theta_0+x_0\\xS\theta_0+yC\theta_0+y_0\\\theta+\theta_0\end{pmatrix},
\end{align*}
with unit element and inverse defined by
\begin{align*}
e:=\begin{pmatrix}0\\0\\0\end{pmatrix}\quad\text{and}\quad
\begin{pmatrix}x\\y\\\theta\end{pmatrix}^{-1}
:=\begin{pmatrix}-xC\theta-yS\theta\\xS\theta-yC\theta_0\\-\theta\end{pmatrix}
\end{align*}
is a group law. Moreover
\begin{align*}
\begin{pmatrix}x_0\\y_0\\\theta_0\end{pmatrix}\begin{pmatrix}u\\v\\x_i\\y_i\\\lambda_i\end{pmatrix}
&:=\begin{pmatrix}u\\v\\x_iC\theta_0-y_iS\theta_0+x_0\\x_iS\theta_0+y_iC\theta_0+y_0\\\lambda_i\end{pmatrix}
\end{align*}
defines a transformation group. Notice the beacon coordinates $(x_i,y_i)$ can be considered as known constant inputs.

Direct computation shows the system is invariant in the sense of definition~\ref{def:invsys}.
Notice $u,v$ are invariant, hence the permanent trajectories are defined by
\begin{align*}
x(t) &=\bar u\sin(\bar u\bar vt+\bar\theta)\\
y(t) &=-\bar u\cos(\bar u\bar vt+\bar\theta)\\
\theta(t) &=\bar u\bar vt,
\end{align*}
where $\bar u,\bar v,\bar\theta$ are arbitrary constants; there consist of arbitrary circles (when $\bar v\neq0$) and lines (when $\bar v=0$).


\subsection{Design of a linearized state controller}
The state tracking error is
\begin{align*}
\begin{pmatrix}\eta_x\\\eta_y\\\eta_\theta\end{pmatrix}
:=\begin{pmatrix}x_r\\y_r\\\theta_r\end{pmatrix}^{-1}\begin{pmatrix}x\\y\\\theta\end{pmatrix}
=\begin{pmatrix}(x-x_r)C\theta_r+(y-y_r)S\theta_r\\-(x-x_r)S\theta_r+(y-y_r)C\theta_r\\\theta-\theta_r\end{pmatrix}
\end{align*}
and satisfies
\begin{align*}
\begin{pmatrix}\dot\eta_x\\\dot\eta_y\\\dot\eta_\theta\end{pmatrix}
&=\begin{pmatrix}C\theta_r & S\theta_r &0\\ -S\theta_r &C\theta_r &0\\ 0&0&1\end{pmatrix}
\begin{pmatrix}\dot x-\dot x_r\\ \dot y-\dot y_r\\ \dot\theta-\dot\theta_r\end{pmatrix}
\\&\quad
+\dot\theta_r\begin{pmatrix}-S\theta_r & C\theta_r &0\\ -C\theta_r &-S\theta_r &0\\ 0&0&0\end{pmatrix}
\begin{pmatrix}x-x_r\\ y-y_r\\0\end{pmatrix}\\
&=\begin{pmatrix}(u_r+\eta_u)C\eta_\theta-u_r+u_rv_r\eta_y\\ 
(u_r+\eta_u)S\eta_\theta-u_rv_r\eta_x\\
(u_r+\eta_u)(v_r+\eta_v)-u_rv_r\end{pmatrix}.
\end{align*}
The linearized error equation around $(\bar\eta_x,\bar\eta_y,\bar\eta_\theta,\bar\eta_u,\bar\eta_v)=(0,0,0,0,0)$ is then
\begin{align*}
\dot\xi_x &=\xi_u+u_rv_r\xi_y\\
\dot\xi_y &=u_r\xi_\theta-u_rv_r\xi_x\\
\dot\xi_\theta &=v_r\xi_u+u_r\xi_v.
\end{align*}

If we choose for instance the state feedback
\begin{align*}
\xi_u &=-u_rv_r\xi_y-\abs{u_r}k_1\xi_x\\
\xi_v &=\sign(u_r)k_1\xi_x+v_r^2\xi_y-k_2\xi_y-\sign(u_r)k_3\xi_\theta
\end{align*}
the resulting closed-loop 
\begin{align*}
\dot\xi_x &=-\abs{u_r}k_1\xi_x\\
\dot\xi_y &=u_r\xi_\theta-u_rv_r\xi_x\\
\dot\xi_\theta &=-u_rk_2\xi_y-\abs{u_r}k_3\xi\theta
\end{align*}
with $k_1,k_2,k_3>0$ is obviously stable for $u_r,v_r$ constant. It is even stable for non constant $u_r,v_r$ provided for $u_r,v_r$ are bounded and $\int_{t_0}^{+\infty}\abs{u_r(t)}dt=+\infty$ for all $t_0\geq0$.

\subsection{Design of an extended observer}
The output observation error is
\begin{align*}
\eps_i &:=\begin{pmatrix}\hat x\\\hat y\\\hat\theta\end{pmatrix}^{-1}
\bigl((\hat x-x_i)^2+(\hat y-y_i)^2\bigr)
-\begin{pmatrix}\hat x\\\hat y\\\hat\theta\end{pmatrix}^{-1}\lambda_i\\
&=(\hat x-x_i)^2+(\hat y-y_i)^2-\lambda_i
\end{align*}
and the ``product''~\eqref{eq:x0xi} is
\begin{align*}
\begin{pmatrix}x_0\\y_0\\\theta_0\end{pmatrix}\xi
&=\begin{pmatrix}C\hat\theta & -S\hat\theta &0\\
S\hat\theta &C\hat\theta &0\\ 0&0&1\end{pmatrix}\xi.
\end{align*}
Every invariant observer then reads
\begin{align*}
\begin{pmatrix}\dot{\hat x}\\ \dot{\hat y}\\ \dot{\hat\theta}\end{pmatrix}
=\begin{pmatrix}uC\hat\theta\\ uS\hat\theta\\ uv\end{pmatrix}
-\begin{pmatrix}C\hat\theta & -S\hat\theta &0\\
S\hat\theta &C\hat\theta &0\\ 0&0&1\end{pmatrix}L
\begin{pmatrix}\eps_1 \\ \vdots\\ \eps_p\end{pmatrix},
\end{align*}
where $L$ is a $3\times p$ matrix possibly depending on $\eps_i$ and 
the invariant quantities
\begin{align*}
\begin{pmatrix}\hat I_u\\\hat I_v\\\hat I_{x_i}\\\hat I_{y_i}\end{pmatrix}
&:=\begin{pmatrix}\hat x\\\hat y\\\hat\theta\end{pmatrix}^{-1}\begin{pmatrix}u\\v\\x_i\\y_i\end{pmatrix}
=\begin{pmatrix}u\\v\\(x_i-\hat x)C\hat\theta+(y_i-\hat y)S\hat\theta\\-(x_i-\hat x)S\hat\theta+(y_i-\hat y)C\hat\theta\end{pmatrix}.
\end{align*}
For reasons that will shortly be apparent, we choose
$$L:=-\frac{1}{2}\LL(\II\II^T)^{-1}\II$$
where $\LL$ is a 3x2 matrix to be defined and $\II$ is the $2\times p$ matrix
\[ \II:=\begin{pmatrix}\hat I_{x_1} &\cdots &\hat I_{x_p}\\ \hat I_{y_1} &\cdots &\hat I_{y_p}\end{pmatrix}. \]
Notice $(\II\II^T)^{-1}$ is always invertible when~$p\geq3$.

The state observation error
\begin{align*}
\begin{pmatrix}\eps_x\\\eps_y\\\eps_\theta\end{pmatrix}
:=\begin{pmatrix}x\\y\\\theta\end{pmatrix}^{-1}\begin{pmatrix}\hat x\\\hat y\\\hat\theta\end{pmatrix}
=\begin{pmatrix}(\hat x-x)C\theta+(\hat y-y)S\theta\\-(\hat x-x)S\theta+(\hat y-y)C\theta\\\hat\theta-\theta\end{pmatrix}
\end{align*}
has for equation
\begin{align*}
\begin{pmatrix}\dot\eps_x\\\dot\eps_y\\\dot\eps_\theta\end{pmatrix}
&=\begin{pmatrix}C\theta & S\theta &0\\ -S\theta &C\theta &0\\ 0&0&1\end{pmatrix}
\left[\begin{pmatrix}\dot{\hat x}-\dot x\\ \dot{\hat y}-\dot y\\ \dot{\hat\theta}-\dot\theta\end{pmatrix}
-\begin{pmatrix}C\hat\theta & -S\hat\theta &0\\ S\hat\theta &C\hat\theta &0\\ 0&0&1\end{pmatrix}L
\begin{pmatrix}\eps_1 \\ \vdots\\ \eps_p\end{pmatrix}\right]
\\&\quad
+\dot\theta\begin{pmatrix}-S\theta & C\theta &0\\ -C\theta &-S\theta &0\\ 0&0&0\end{pmatrix}
\begin{pmatrix}\hat x-x\\\hat y-y\\0\end{pmatrix}\\
&=\begin{pmatrix}u(C\eps_\theta-1)+uv\eps_y\\ -uS\eps_\theta-uv\eps_x\\0\end{pmatrix}
+\begin{pmatrix}C\eps_\theta & -S\eps_\theta &0\\ S\eps_\theta &C\eps_\theta &0\\ 0&0&1\end{pmatrix}L
\begin{pmatrix}\eps_1 \\ \vdots\\ \eps_p\end{pmatrix},
\end{align*}
with
\begin{align*}
\eps_i &=(\hat x-x_i)^2+(\hat y-y_i)^2-(x-x_i)^2-(y-y_i)^2\\
&=-2\begin{pmatrix}x_i-\hat x &y_i-\hat y\end{pmatrix}\begin{pmatrix}\hat x-x\\\hat y-y\end{pmatrix}-(\hat x-x)^2-(\hat y-y)^2\\
&=-2\begin{pmatrix}\hat I_{x_i} &\hat I_{y_i}\end{pmatrix}
\begin{pmatrix}C\hat\theta &S\hat\theta\\-S\hat\theta &C\hat\theta\end{pmatrix}
\begin{pmatrix}C\theta &-S\theta\\S\theta &C\theta\end{pmatrix}
\begin{pmatrix}\eps_x\\\eps_y\end{pmatrix}-\eps_x^2-\eps_y^2\\
&=-2\begin{pmatrix}\hat I_{x_i} &\hat I_{y_i}\end{pmatrix}
\begin{pmatrix}C\eps_\theta &-S\eps_\theta\\S\eps_\theta &C\eps_\theta\end{pmatrix}
\begin{pmatrix}\eps_x\\\eps_y\end{pmatrix}-\eps_x^2-\eps_y^2.
\end{align*}
Linearizing around the equilibrium point $(\bar\eps_x,\bar\eps_y,\bar\eps_\theta):=(0,0,0)$, we have
\begin{align*}
\begin{pmatrix}\dot e_x\\ \dot e_y\\ \dot e_\theta\end{pmatrix}
&=\begin{pmatrix}uve_y\\ -ue_\theta-uve_x\\0\end{pmatrix}
-L\begin{pmatrix}e_1 \\ \vdots\\ e_p\end{pmatrix},
\end{align*}
with
\begin{align*}
e_i &=-2\begin{pmatrix}\hat I_{x_i} &\hat I_{y_i}\end{pmatrix}
\begin{pmatrix}e_x\\e_y\end{pmatrix}.
\end{align*}
Hence
\begin{align*}
\begin{pmatrix}e_1 \\ \vdots\\ e_p\end{pmatrix}&=-2\II^T\begin{pmatrix}e_x\\ e_y\end{pmatrix}
\end{align*}
so we eventually have
\begin{align*}
\begin{pmatrix}\dot e_x\\ \dot e_y\\ \dot e_\theta\end{pmatrix}
&=\begin{pmatrix}uve_y\\ -ue_\theta-uve_x\\0\end{pmatrix}
-\LL\begin{pmatrix}e_x\\ e_y\end{pmatrix}.
\end{align*}
If we choose
$$\LL:=\begin{pmatrix}\abs{u}l_1 &uv\\ -uv &\abs{u}l_2\\ 0 &-ul_3\end{pmatrix}$$
with $l_1,l_2,l_3>0$, the linearized error system reads
\begin{align*}
\dot e_x &=-\abs{u}l_1e_x\\
\dot e_y &=-ue_\theta-\abs{u}l_2e_y\\ 
\dot e_\theta &=-ul_3e_y
\end{align*}
is obviously stable provided $\int_{t_0}^{+\infty}\abs{u_r(t)}dt=+\infty$ for all $t_0\geq0$. This means the designed observer is convergent locally around \emph{every} trajectory.

\subsection{The control law}
Combining as in section~\ref{sec:extobsLie} the previously designed linearized controller and extended observer
we end up with the control law
\begin{align*}
\begin{pmatrix}u\\ v\end{pmatrix}
&=\begin{pmatrix}u_r-u_rv_r\hat\eta_y-\abs{u_r}k_1\hat\eta_x\\
v_r\sign(u_r)k_1\hat\eta_x+v_r^2\xi_y-k_2\hat\eta_y-\sign(u_r)k_3\hat\eta_\theta\end{pmatrix}\\
\begin{pmatrix}\dot{\hat x}\\ \dot{\hat y}\\ \dot{\hat\theta}\end{pmatrix}
&=\begin{pmatrix}uC\hat\theta\\ uS\hat\theta\\ uv\end{pmatrix}
+\frac{1}{2}\begin{pmatrix}C\hat\theta & -S\hat\theta &0\\
S\hat\theta &C\hat\theta &0\\ 0&0&1\end{pmatrix}\LL(\II\II^T)^{-1}\II
\begin{pmatrix}\eps_1 \\ \vdots\\ \eps_p\end{pmatrix},
\end{align*}
where
\begin{align*}
\begin{pmatrix}\hat\eta_x\\\hat\eta_y\\\hat\eta_\theta\end{pmatrix}
:=\begin{pmatrix}x_r\\y_r\\\theta_r\end{pmatrix}^{-1}\begin{pmatrix}\hat x\\\hat y\\\theta\end{pmatrix}
=\begin{pmatrix}(\hat x-x_r)C\theta_r+(\hat y-y_r)S\theta_r\\-(\hat x-x_r)S\theta_r+(\hat y-y_r)C\theta_r\\\hat\theta-\theta_r\end{pmatrix}.
\end{align*}
Notice the map $\kappa$ chosen here is simply the identity. This control law stabilizes the system around any permanent trajectory.
\section{Fully-actuated mechanical systems on Lie groups}\label{bullo:sec}

In this section, we consider a simple mechanical system on a Lie group whose motion is described by the so-called Euler-Poincar\'e equations. We focus on fully-actuated systems arising in control theory as described by e.g. \cite{bullo-murray-auto99}
 \begin{align}
\dot x&=f(x,\xi)\label{euler-poincare:eq1}\\
\dot \xi&=A(\xi)+I^{-1}(F(x,\xi)+u)\label{euler-poincare:eq2}
\end{align}
where $x\in G$ is the (generalized) position, with $G$ a Lie group (the configuration space), $\xi\in T_xG$ is the (generalized) velocity, $A$ is a bilinear function of its argument, $I^{-1}(F(x,\xi)+u)$ is the resultant force acting on the system, and $u\in T_xG$ denotes the control. Moreover, the system  \eqref{euler-poincare:eq1} is invariant to the following $G$-group action : for any $x_0\in G$ the action $x_0x$ is the left multiplication on $G$, and  $x_0\xi=\xi$.

\subsection{Considered observation problem}
We assume that $\xi$ is known or measured,  and dropping the second equation \eqref{euler-poincare:eq2} (as $\xi$ does not need to be estimated), we focus on the following observation problem:
 \begin{equation}\label{red:eq}\begin{aligned}
\dot x&=f(x,\xi)\\
y&=h(x)
\end{aligned}\end{equation}
where the output satisfies $x_0h(x)=h(x_0x)$. The theory of Section 2 applies to this problem, and the permanent trajectories are generated by constant velocities, i.e.  $I_r=x_r^{-1}\xi_r=\xi_r\equiv cst$. Such trajectories are very natural for mechanical systems and admit a geometrical interpretation (they are generated by one-parameter subgroups of $G$). They constitute interesting motion primitives  that can be concatenated to yield a very large class of trajectories (think of straight lines and coordinated turn in avionics which are permanent trajectories on $SE(3)$ \cite{arxiv-08}). In the sequel we will prove that a local separation principle holds around those trajectories.

\subsection{Control problem around permanent trajectories}

Consider the reference trajectory ($x_r(t),\xi_r$) generated by a time-invariant $\xi_r$. Let $u_r$ be the corresponding control.  Let the invariant tracking error be $(\eta_x,\eta_\xi)=(x_r^{-1}x,\xi-\xi_r)$ as in \cite{bullo-murray-auto99}. As the system \eqref{euler-poincare:eq1} is invariant, applying \eqref{eq:etaxLie}  (where $u$ is replaced by $\xi$) we have 
$$
\dotex \eta_x=\Upsilon(\eta_x,\eta_\xi,I_r)
$$
Moreover we have, up to second order terms in $(\eta_x,\eta_\xi)$
\begin{align*}
\dotex \eta_\xi&=A(I_r+\eta_\xi)-A(I_r)+I^{-1}(\delta_1{F}({x_r,I_r})\eta_x\\&\quad+\delta_2{F}({x_r,I_r})\eta_\xi+\eta_u)
\end{align*}We have thus the following result:
\begin{lem}\label{lem:autonomous}
If $\delta_1{F}({x,\xi})$ and $\delta_2{F}({x,\xi})$ do not depend on $x$, the linearized tracking error around permanent trajectories generated by $\xi_r\equiv cst$ is \emph{time-invariant}.
\end{lem}
It implies that under assumptions of Lemma \ref{lem:autonomous} (which are usually satisfied in practice), a local separation principle holds around permanent trajectories for the control problem \eqref{euler-poincare:eq1}-\eqref{euler-poincare:eq2} combined with an invariant observer for the subsystem \eqref{red:eq}, and where $\xi$ is known or measured.

We believe that the problem addressed in this section is relevant to applications, and examples are left for future research. Note that  \cite{maithripala2005}  proved a separation principle for simple mechanical systems on Lie groups where the position is measured, and a {velocity observer}  must be designed. This is a different problem and a different approach is developed.  The drawback of this approach is that it relies on many restrictive assumptions. In particular:
the group must be compact, and the observer must be exponentially convergent around $\emph{any}$ trajectory.

\section{Conclusion}

In this paper a local separation principle around a large set of trajectories for non-linear invariant systems on Lie groups was proved. We also proved that the results extend partially to the control of simple mechanical systems on Lie groups. In future research we plan to explore examples of mechanical systems for which those results apply. 


\end{document}